\documentstyle[amscd]{amsart}

\newenvironment{proof}{\noindent{\bf Proof.}}{\qed\par}
\newenvironment{proofof}[1]{\noindent{\bf Proof of #1.}}{\qed\par}

\theoremstyle{plain}
\newtheorem{thm}{Theorem}[subsection]
\newtheorem{pro}[thm]{Proposition}
\newtheorem{cor}[thm]{Corollary}
\newtheorem{lem}[thm]{Lemma}

\theoremstyle{definition}
\newtheorem{defn}[thm]{Definition}

\theoremstyle{remark}
\newtheorem{remark}[thm]{Remark}

\newcommand{\oper}{\operatorname}

\newcommand{\DD}[2]{{\frak D}_{#1,#2}}
\newcommand{\fL}[2]{{\frak L}_{#1,#2}}
\newcommand{\MM}[1]{{\cal Q}_{#1}}
\newcommand{\MMK}[1]{{\cal Q}_{#1}^K}
\newcommand{\MMD}[1]{{\cal Q}_{#1}^D}
\newcommand{\MML}[1]{{\cal Q}_{#1}^L}
\newcommand{\BL}{{\mbox{\bf L}}}
\newcommand{\bN}{{\mbox{\bf N}}}
\newcommand{\bB}{{\mbox{\bf B}}}
\newcommand{\LL}{{\cal L}}
\newcommand{\CA}{{\cal A}}
\newcommand{\CV}{{\cal V}}
\newcommand{\cV}{{\cal V}}
\newcommand{\cX}{{\cal X}}
\newcommand{\cY}{{\cal Y}}
\newcommand{\CH}{{\cal H}}

\newcommand{\CF}{{\cal F}}
\newcommand{\NN}{{\Bbb N}}
\newcommand{\CC}{{\Bbb C}}
\newcommand{\AAA}{{\Bbb A}}
\newcommand{\KK}{{\frak K}}
\newcommand{\fk}{{\frak k}}
\newcommand{\fA}{{\frak A}}
\newcommand{\MMM}{{\frak M}}
\newcommand{\SS}{{\frak S}}
\newcommand{\TT}{{\frak T}}
\newcommand{\ZZ}{{\Bbb Z}}
\newcommand{\QQ}{{\Bbb Q}}
\newcommand{\UQ}{\underline{\Bbb Q}}
\newcommand{\FF}{{\cal B}}
\newcommand{\PP}{{\Bbb P}}

\newcommand{\EE}{{\cal E}}
\newcommand{\OO}{{\cal O}}
\newcommand{\cK}{{\cal K}}
\newcommand{\mm}{{\frak m}}
\newcommand{\om}{\omega}
\newcommand{\ka}{\kappa}
\newcommand{\eps}{\varepsilon}
\newcommand{\HH}{\oper{H}\nolimits}

\newcommand{\vphi}{\varphi}
\newcommand{\ti}{\tilde}
\newcommand{\rnk}{l}
\newcommand{\leng}{\ell}
\newcommand{\lra}{\longrightarrow}
\newcommand{\ds}{\displaystyle}
\newcommand{\suml}{\sum\limits}

\newcommand{\refeq}[1]{{\rm (\ref{#1})}}

\newcommand{\rank}{\oper{rank}}
\newcommand{\codim}{\oper{codim}}
\newcommand{\deff}{\oper{def}}
\newcommand{\Ker}{\oper{Ker}}
\newcommand{\Hom}{\oper{Hom}}
\newcommand{\Ext}{\oper{Ext}\nolimits}

\newcommand{\ede}{\stackrel{\text{\rm def}}{=}}
\newcommand{\coin}[3]{\nabla_{#1}\left(#2,#3\right)}
\newcommand{\pq}[1]{{1\le q\le p\le #1}}

\title{Laumon's resolution of Drinfeld's compactification is small}
\author{Alexander Kuznetsov}
\address{Independent University of Moscow, 13-26, Vereskovaja St.,
Moscow 129329 Russia}
\email{sasha@@ium.ips.ras.ru}
\thanks{I was partially supported by AMS and CRDF-$2819$ grants.}

\begin{document}

\maketitle

Let $C$ be a smooth projective curve of genus $0$.
Let $\FF$ be the variety of complete flags in an
$n$-dimensional vector space $V$. Given an
$(n-1)$-tuple $\alpha$ of positive integers one
can consider the space $\MM\alpha$ of algebraic
maps of degree $\alpha$ from $C$ to $\FF$. This space has drawn much
attention recently in connection with Quantum
Cohomology (see e.g. \cite{Givental}, \cite{Kontsevich}). The
space $\MM\alpha$ is smooth but not compact
(see e.g. \cite{Kontsevich}). The problem of compactification
of $\MM\alpha$ proved very important. One
compactification $\MMK\alpha$ was constructed in
loc.\ cit.\ (the space of {\em stable maps}).
Another compactification $\MML\alpha$ (the space
of {\em quasiflags}), was constructed in \cite{Laumon}.
However, historically the first and most economical
compactification $\MMD\alpha$ (the space of {\em
quasimaps}) was constructed by Drinfeld (early 80-s,
unpublished). The latter compactification is singular,
while the former ones are smooth. Drinfeld has
conjectured that the natural map $\pi:\MML\alpha\to\MMD\alpha$
is a small resolution of singularities. In the present
note we prove this conjecture after the necessary
recollections. In fact, the proof gives some additional
information about the fibers of $\pi$. It appears that
every fiber has a cell decomposition, i.e. roughly
speaking, is a disjoint union of affine spaces. This
permits to compute not only the stalks of $IC$ sheaf on
$\MMD\alpha$ but, moreover, the Hodge structure in
these stalks. Namely, the Hodge structure is a pure Tate one,
and the generating function for the $IC$ stalks is just
the Lusztig's $q$-analogue of Kostant's partition function
(see \cite{Lusztig}).

In conclusion, let us mention that the Drinfeld compactifications are
defined for the space of maps into flag manifolds of arbitrary
semisimple group, and we plan to study their small resolutions in
a paper to follow.

I am very grateful to M.~Finkelberg for very useful discussions
and permanent attention.
I would like to thank D.~Gaitsgory who has explained to me the
statement of Drinfeld's conjecture.

\section{The space of maps into flag variety}

\subsection{Notations}

Let $G$ be a complex semisimple simply-connected Lie group,
$H\subset B$ its Cartan and Borel subgroups,
$N$ the unipotent radical of $B$,
$Y$ the lattice of coroots of $G$ (with respect to $H$),
$\rnk$ the rank of $Y$,
$I=\{i_1,i_2,\dots i_\rnk\}$ the set of simple coroots,
$R^+$ the set of positive coroots,
$X$ the lattice of weights,
$X^+$ the cone of dominant weights,
$\Omega=\{\om_1,\om_2\dots\om_\rnk\}$ the set of fundamental weights
($\langle\om_k,i_l\rangle=\delta_{kl}$),
$\FF=G/B$ the flag variety and
$C$ a smooth projective curve of genus ~$0$.
Recall that there are canonical isomorphisms
$$
\HH_2(\FF,\ZZ)\cong Y\qquad\HH^2(\FF,\ZZ)\cong X.
$$
For $\lambda\in X$ let $\BL_\lambda$ denote the corresponding
$G$-equivariant line bundle on $\FF$.

The map $\vphi:C\to\FF$ has degree $\alpha\in\NN[I]\subset Y$
if the following equivalent conditions hold:
\begin{enumerate}
\item $\vphi_*([C])=\alpha$;
\item for any $\lambda\in X\quad\text{we have}\quad
\deg(\vphi^*\BL_\lambda)=\langle\lambda,\alpha\rangle$.
\end{enumerate}

We denote by $\MM\alpha$ the space of algebraic maps from $C$
to $\FF$ of degree $\alpha$.
It is known that $\MM\alpha$ is smooth variety and
$\dim\MM\alpha=2|\alpha|+\dim\FF$. In this paper we compare
two natural compactifications of the space $\MM\alpha$, which
we presently describe.

\subsection{Drinfeld's compactification}

The Pl\"ucker embedding of the flag variety $\FF$ gives rise to the
following interpretation of $\MM\alpha$.

For any irreducible representation $V_\lambda$ ($\lambda\in X^+$) of $G$
we consider the trivial vector bundle $\CV_\lambda=V_\lambda\otimes\OO_C$
over $C$.

For any $G$-morphism $\psi:\ V_\lambda\otimes V_\mu\lra V_\nu$ we denote
by the same letter the induced morphism $\psi:\ \CV_\lambda\otimes \CV_\mu
\lra \CV_\nu$.

Then $\MM\alpha$ is the space of collections of line subbundles
$\LL_\lambda\subset\CV_\lambda,\ \lambda\in X^+$ such that:

a) $\deg\LL_\lambda=-\langle\lambda,\alpha\rangle$;

b) For any nonzero $G$-morphism $\psi:\ V_\lambda\otimes V_\mu\lra V_\nu$
such that $\nu=\lambda+\mu$ we have $\psi(\LL_\lambda\otimes\LL_\mu)=
\LL_\nu$;

c) For any $G$-morphism $\psi:\ V_\lambda\otimes V_\mu\lra V_\nu$
such that $\nu<\lambda+\mu$ we have $\psi(\LL_\lambda\otimes\LL_\mu)=0$.

\begin{remark}\label{rem1}
Certainly, the property b) guarantees that in order to specify such a
collection it suffices to give $\LL_{\omega_k}$ for the set $\Omega$ of
fundamental weights.
\end{remark}

If we replace the curve $C$ by a point, we get the Pl\"ucker description of
the flag variety $\FF$ as the space of collections of lines
$L_\lambda\subset V_\lambda$ satisfying conditions of type (b) and (c)
(thus $\FF$ is embedded into $\prod\limits_{\lambda\in X^+}\PP(V_\lambda)$).
Here, a Borel subgroup $\bB$ in $\FF$ corresponds to
a system  of lines $(L_\lambda,\lambda\in X^+)$ if lines are the fixed
points of the unipotent radical of $\bB$, $L_\lambda=(V_\lambda)^\bN$,
or equivalently, if $\bN$ is the common stabilizer for all lines
$\bN=\bigcap\limits_{\lambda\in X^+}G_{L_\lambda}$.

The following definition in case $G=SL_2$ appeared in \cite{Drinfeld}.

\begin{defn}[V.Drinfeld]\label{MMD}
The space $\MMD\alpha$ of quasimaps of degree $\alpha$
from $C$ to $\FF$ is the space of
collections of invertible subsheaves
$\LL_\lambda\subset\CV_\lambda,\ \lambda\in X^+$ such that:

a) $\deg\LL_\lambda=-\langle\lambda,\alpha\rangle$;

b) For any nonzero $G$-morphism $\psi:\ V_\lambda\otimes V_\mu\lra V_\nu$
such that $\nu=\lambda+\mu$ we have $\psi(\LL_\lambda\otimes\LL_\mu)=
\LL_\nu$;

c) For any $G$-morphism $\psi:\ V_\lambda\otimes V_\mu\lra V_\nu$
such that $\nu<\lambda+\mu$ we have $\psi(\LL_\lambda\otimes\LL_\mu)=0$.
\end{defn}

\begin{remark}
Here is another version of the Definition, also due to V.Drinfeld.
The principal affine space $\CA=G/N$ is an $H$-torsor over $\FF$.
We consider its affine closure $\fA$,
that is, the spectrum of the ring of functions on $\CA$.
The action of $H$ extends to $\fA$ but it is not free anymore. Consider the
quotient stack $\ti\FF=\fA/H$. The flag variety $\FF$ is an
open substack in $\ti\FF$. A map
$\ti\phi:C\to\ti\FF$ is nothing else than an $H$-torsor $\Phi$ over $C$ along
with an $H$-equivariant morphism $f:\Phi\to\fA$. The degree of this map
is defined as follows.

Let $\chi_\lambda:H\to\CC^*$ be the character of $H$ corresponding
to a weight $\lambda\in X$. Let $H_\lambda\subset H$ be the kernel
of the morphism $\chi_\lambda$. Consider the induced $\CC^*$-torsor
$\Phi_\lambda=\Phi/H_\lambda$ over $C$. The map $\ti\phi$ has
degree $\alpha\in\NN[I]$ if
$$
\text{for any }\lambda\in X\quad\text{we have}\quad
\deg(\Phi_\lambda)=\langle\lambda,\alpha\rangle.
$$
\end{remark}
\begin{defn}\label{MMD1}
The space $\MMD\alpha$ is the space of maps $\ti\phi:C\to\ti\FF$
of degree $\alpha$ such that the generic point of $C$ maps into
$\FF\subset\ti\FF$.
\end{defn}
The equivalence of \ref{MMD} and \ref{MMD1} follows immediately from
the Pl\"ucker embedding of $\fA$ into
$\prod\limits_{\lambda\in X^+}V_\lambda$.

\begin{pro}
$\MMD\alpha$ is a projective variety.
\end{pro}
\begin{proof}
The space $\MMD\alpha$ is naturally embedded into the space
$$
\prod_{k=1}^\rnk
\PP(\Hom(\OO_C(-\langle\om_k,\alpha\rangle),\CV_{\om_k}))
$$
and is closed in it.
\end{proof}

\subsection{The stratification of the Drinfeld's compactification}

In this subsection we will introduce the stratification of the space
of quasimaps.

{\bf Configurations of $I$-colored divisors.}

Let us fix $\alpha\in\NN[I]\subset Y,\ \alpha=\suml_{k=1}^N a_ki_k$.
Consider the configuration space $C^\alpha$ of colored
effective divisors of multidegree $\alpha$ (the set of colors is $I$).
The dimension of $C^\alpha$ is equal to the length
$|\alpha|=\suml_{k=1}^N a_k$.

Multisubsets of a set $S$ are defined as elements of some
symmetric power $S^{(m)}$ and we denote the image of
$(s_1,\dots,s_m)\in S^m$ by $\{\{s_1,\dots,s_m\}\}$.
We denote by $\Gamma(\alpha)$ the set of all partitions of $\alpha$,
i.e.\ multisubsets $\Gamma=\{\{\gamma_1,\dots,\gamma_m\}\}$ of $\NN[I]$
with $\suml_{r=1}^m\gamma_r=\alpha$, $\gamma_r>0$.

For $\Gamma\in\Gamma(\alpha)$ the corresponding stratum $C^\alpha_\Gamma$
is defined as follows. It is formed by configurations which can be
subdivided into $m$ groups of points, the $r$-th group containing $\gamma_r$
points; all the points in one group equal to each other, the different
groups being disjoint. For example, the main diagonal in $C^\alpha$
is the closed stratum given by partition $\alpha=\alpha$, while the complement
to all diagonals in $C^\alpha$ is the open stratum given by partition
$$
\alpha=\suml_{k=1}^N(\underbrace{i_k+i_k+\ldots+i_k}_{a_k\oper{ times}})
$$
Evidently, $C^\alpha=\bigsqcup\limits_{\Gamma\in\Gamma(\alpha)}C^\alpha_\Gamma$.

{\bf Normalization and defect of subsheaves.}\nopagebreak

Let $F$ be a vector bundle on the curve $C$ and let $E$ be a subsheaf in $F$.
Let $F/E=T(E)\oplus L$ be the decomposition of the quotient sheaf $F/E$
into the sum of its torsion subsheaf and a locally free sheaf, and let
$\ti E=\Ker(F\to L)$ be the kernel of the natural map $F\to L$.
Then $\ti E$ is a vector subbundle in $F$ which contains $E$ and has
the following universal property:
$$
\text{for any subbundle $\EE'\subset F$ if $\EE'$ contains $E$
then $\EE'$ contains also $\ti E$}.
$$
Moreover, $\rank\ti E=\rank E$, $\ti E/E\cong T(E)$ and
$c_1(\ti E)=c_1(E)+\leng(T(E))$ (for any torsion sheaf on $C$
we denote by $\leng(T)$ its length).
\begin{defn}
We will call $\ti E$ the {\em normalization}
of $E$ in $F$ and $T(E)$ the {\em defect} of $E$.
\end{defn}
\begin{remark}\label{detnorm}
If $\ti E$ is the normalization of $E$ in $F$ then $\Lambda^k(\ti E)$ is the
normalization of $\Lambda^kE$ in $\Lambda^k F$.
\end{remark}

For any $x\in C$ and torsion sheaf $T$ on $C$ we will denote by
$\leng_x(T)$ the length of the localization of $T$ in the point $x$.
\begin{defn}
For any quasimap
$\vphi=(\LL_\lambda\subset\CV_\lambda)_{\lambda\in X^+}\in\MMD\alpha$
we define the {\em normalization} of $\vphi$ as follows:
$$
\ti\vphi=(\ti\LL_\lambda\subset\CV_\lambda),
$$
and the {\em defect} of $\vphi$ as follows:
$$
\deff(\vphi)=(T(\LL_\lambda))
$$
(the defect of $\vphi$ is a collection of torsion sheaves).
\end{defn}
\begin{pro}\label{normalization}
For any $\vphi\in\MMD\alpha$ there exists $\beta\le\alpha\in\NN[I]$,
partition $\Gamma=(\gamma_1,\dots,\gamma_m)\in\Gamma(\alpha-\beta)$ and
a divisor $D=\suml_{r=1}^m\gamma_rx_r\in C^{\alpha-\beta}_\Gamma$ such that
$$
\ti\vphi\in\MM\beta
\qquad
\leng_x(\deff(\vphi)_\lambda)=\begin{cases}
\ds\langle\lambda,\gamma_r\rangle,&\text{if $x=x_r$}\\
\ds0,&\text{otherwise}\end{cases}
$$
\end{pro}
\begin{proof}
Clear.
\end{proof}

\begin{defn}
The pair $(\beta,\Gamma)$ will be called the {\em type of degeneration}
of $\vphi$. We denote by $\DD\beta\Gamma$ the subspace of $\MMD\alpha$
consisting of all quasimaps $\vphi$ with the given type of degeneration.
\end{defn}
\begin{remark}
Note that $\DD0\emptyset=\MM\alpha$.
\end{remark}

We have
\begin{equation}
\MMD\alpha=\bigsqcup\begin{Sb}
\beta\le\alpha\\
\Gamma\in\Gamma(\alpha-\beta)
\end{Sb}\DD\beta\Gamma
\end{equation}

The map $d_{\beta,\Gamma}:
\DD\beta\Gamma\to\MM\beta\times C^{\beta-\alpha}_\Gamma$ which sends
$\vphi$ to $(\ti\vphi,D)$ (see \ref{normalization}) is an isomorphism.
The inverse map $\sigma_{\beta,\Gamma}$ can be constructed as follows.
Let $\vphi=(\LL_\lambda)\in\MM\beta$. Then
$$
\sigma_{\beta,\Gamma}(\vphi,D)\ede(\LL'_\lambda)\qquad
\LL'_\lambda\ede
\bigcap_{r=1}^m\mm_{x_r}^{\langle\lambda,\gamma_r\rangle}\cdot\LL_\lambda,
$$
where $\mm_x$ denotes the sheaf of ideals of the point $x\in C$.

\subsection{Laumon's compactification}

Let $V$ be an $n$-dimensional vector space.\linebreak
>From now on we will assume that $G=SL(V)$ (in this case
certainly $\rnk=n-1$). In this case there is the Grassmann
embedding of the flag variety, namely
$$
\FF=\{(U_1,U_2,\dots,U_{n-1})\in
G_1(V)\times G_2(V)\times\dots\times G_{n-1}(V)\ |\
U_1\subset U_2\subset\dots\subset U_{n-1}\},
$$
where $G_k(V)$ is the Grassmann variety of $k$-dimensional subspaces in $V$.
This embedding gives rise to another interpretation of $\MM\alpha$.

We will denote by $\CV$ the trivial vector bundle $V\otimes\OO_C$ over $C$.
Let $\alpha=\suml_{k=1}^{n-1}a_ki_k$, where $i_k$ is the simple coroot dual
to the highest weight $\om_k$ of representation $G$ in $\Lambda^kV$.

Then $\MM\alpha$ is the space of complete flags of vector subbundles
$$
0\subset\EE_1\subset\EE_2\subset\dots\subset\EE_{n-1}\subset\CV
\quad\text{ such that }\quad c_1(\EE_k)=-\langle\om_k,\alpha\rangle=-a_k.
$$
\begin{defn}[{Laumon, \cite[4.2]{Laumon}}]
The space $\MML\alpha$ of {\em quasiflags} of degree $\alpha$ is the space
of complete flags of locally free subsheaves
$$
0\subset E_1\subset E_2\subset\dots\subset E_{n-1}\subset\CV
\quad\text{ such that }\quad c_1(E_k)=-\langle\om_k,\alpha\rangle=-a_k.
$$
\end{defn}

It is known that $\MML\alpha$ is a smooth projective variety
of dimension $2|\alpha|+\dim\FF$ (see loc.\ cit., Lemma 4.2.3).

\subsection{The stratification of the Laumon's compactification}

There is a stratification of the space $\MML\alpha$ similar to
the above stratification of $\MMD\alpha$.
\begin{defn}
For any quasiflag $E_\bullet=(E_1,\dots,E_{n-1})$ we define its
{\em normalization} as
$$
\ti E_\bullet=(\ti E_1,\dots,\ti E_{n-1}),\text{ where }\ti E_k\text{ is the
normalization of $E_k$ in $\CV$}
$$
and {\em defect}
$$
\deff(E_\bullet)=(\ti E_1/E_1,\dots,\ti E_{n-1}/E_{n-1})
$$
Thus, the defect of $E_\bullet$ is a collection of torsion sheaves.
\end{defn}

\begin{pro}\label{Lnormalization}
For any $E_\bullet\in\MML\alpha$ there exist $\beta\le\alpha\in\NN[I]$,
partition $\Gamma=(\gamma_1,\dots,\gamma_m)\in\Gamma(\alpha-\beta)$ and
a divisor $D=\suml_{r=1}^m\gamma_rx_r\in C^{\alpha-\beta}_\Gamma$ such that
$$
\ti E_\bullet\in\MM\beta,
\qquad
\leng_x(\deff(E_k))=\begin{cases}
\langle\om_k,\gamma_r\rangle,&\text{if $x=x_r$}\\
0,&\text{otherwise}\end{cases}
$$
\end{pro}

\begin{defn}
The pair $(\beta,\Gamma)$ will be called the {\em type of degeneration}
of $E_\bullet$. We denote by $\fL\beta\Gamma$ the subspace in $\MML\alpha$
consisting of all quasiflags $E_\bullet$ with the given type of degeneration.
\end{defn}

\begin{remark}
Note that $\fL0\emptyset=\MM\alpha$.
\end{remark}
We have
\begin{equation}
\MML\alpha=\bigsqcup\begin{Sb}
\beta\le\alpha\\
\Gamma\in\Gamma(\alpha-\beta)
\end{Sb}\fL\beta\Gamma
\end{equation}

\subsection{The map from $\MML\alpha$ to $\MMD\alpha$}

Consider the map $\pi:\MML\alpha\to\MMD\alpha$ which
sends a quasiflag of degree $\alpha$ $E_\bullet\in\MML\alpha$
to a quasimap given by the collection $(\LL_{\om_k})_{k=1}^{n-1}$
(see Remark \ref{rem1}) where
$\LL_{\om_k}=\Lambda^k E_k\subset\Lambda^k\CV=\CV_{\om_k}$.

\begin{pro}
Let $E_\bullet$ be a quasiflag of degree $\alpha$ and let $(\beta,\Gamma)$
be its type of degeneration. Then $\pi(E_\bullet)$ is a quasimap of
degree $\alpha$ and its type of degeneration is $(\beta,\Gamma)$.
\end{pro}
\begin{proof}
Obviously we have
$\deg\LL_{\om_k}=\deg\Lambda^k E_k=c_1(E_k)=-\langle\om_k,\alpha\rangle$
which means that $\pi(E_\bullet)\in\MMD\alpha$.
According to the Remark \ref{detnorm},
$\ti\LL_{\om_k}=\Lambda^k\ti E_k$ (i.e. $(\LL_{\om_k})\in\MM\beta$), hence
\begin{equation}\label{lx}
\leng_x(\ti\LL_{\om_k}/\LL_{\om_k})=\leng_x(\ti E_k/E_k).
\end{equation}
This proves the Proposition.
\end{proof}
\begin{remark}
Note that \refeq{lx} implies that
$\pi$ preserves not only $\beta$ and $\Gamma$ but also $D$ (see
\ref{normalization}, \ref{Lnormalization}).
\end{remark}

Recall that a proper birational map $f:\cX\to\cY$ is called {\em small}
if the following condition holds: let $\cY_m$ be the set of all points
$y\in\cY$ such that $\dim f^{-1}(y)\ge m$. Then for $m>0$ we have
\begin{equation}\label{small}
\codim\cY_m>2m.
\end{equation}

{\bf Main Theorem. } The map $\pi$ is a small resolution of singularities.

\section{The fibers of $\pi$}

\subsection{}
We fix $\EE_\bullet\in\MM\beta$, a partition $\Gamma\in\Gamma(\alpha-\beta)$,
and a divisor $D\in C_\Gamma^{\alpha-\beta}$. Then
$(\EE_\bullet,D)\in\DD\beta\Gamma$. We define $F(\EE_\bullet,D)$
as $\pi^{-1}(\EE_\bullet,D)$.

Let $D=\sum_{r=1}^m\gamma_rx_r$. We define the space $\CF(\EE_\bullet,D)$
of commutative diagrams
$$
\begin{CD}
\EE_1   @>>>    \EE_2   @>>>    \dots   @>>>    \EE_{n-1}   \\
@V\eps_1VV  @V\eps_2VV  @.      @V\eps_{n-1}VV  \\
T_1  @>\tau_1>> T_2  @>\tau_2>> \dots @>\tau_{n-2}>>T_{n-1}
\end{CD}
$$
such that

a) $\eps_k$ is surjective,

b) $T_k$ is torsion,

c) $\ds \leng_x(T_k)=\begin{cases}
\ds\langle\om_k,\gamma_r\rangle,&\text{if $\ds x=x_r$}\\
\ds0,&\text{otherwise}\end{cases}$
\begin{lem}\label{torsion}
We have an isomorphism
$$
F(\EE_\bullet,D)\cong\CF(\EE_\bullet,D).
$$
\end{lem}

\begin{proof}
If $E_\bullet\in F(\EE_\bullet,D)$ then by the \ref{Lnormalization}
the collection $(T_1,\dots,T_{n-1})=\deff(E_\bullet)$ satisfies
the above conditions.

Vice versa, if the collection $(T_1,\dots,T_k)$ satisfies the
above conditions, then consider
$$
E_k=\Ker(\EE_k @>\eps_k>> T_k).
$$
Since the square
$$
\begin{CD}
\EE_k       @>\eps_k>>  T_k     \\
@VVV                @V\tau_kVV  \\
\EE_{k+1}   @>\eps_{k+1}>>  T_{k+1}
\end{CD}
$$
commutes, we can extend it to the commutative diagram
$$
\begin{CD}
0   @>>>    E_k @>>>    \EE_k    @>\eps_k>> T_k @>>>    0 \\
@.      @VVV        @VVV            @V\tau_kVV  @.\\
0   @>>>    E_{k+1} @>>>    \EE_{k+1}@>\eps_{k+1}>> T_{k+1} @>>>    0
\end{CD}
$$
The induced morphism $E_k\to E_{k+1}$ is injective because
$\EE_k\to\EE_{k+1}$ is, and
\begin{multline*}
\qquad c_1(E_k)=c_1(\EE_k)-\leng(T_k)=
-\langle\om_k,\beta\rangle-\sum_{x\in C}\leng_x(T_k)=\\=
-\langle\om_k,\beta\rangle-\sum_{r=1}^m\langle\om_k,\gamma_r\rangle=
-\langle\om_k,\beta+(\alpha-\beta)\rangle=-\langle\om_k,\alpha\rangle\qquad
\end{multline*}
This means that $E_\bullet\in F(\EE_\bullet,D)$.
\end{proof}

\begin{pro}\label{gfibre}
If $D=\suml_{r=1}^m\gamma_rx_r$ is a decomposition
into disjoint divisors then
\begin{equation}
F(\EE_\bullet,D)\cong\prod_{r=1}^mF(\EE_\bullet,\gamma_rx_r).
\end{equation}
\end{pro}
\begin{proof}
Recall that if $T$ is a torsion sheaf on the curve $C$ then
$$
T=\bigoplus_{x\in C}T_x,
$$
where $T_x$ is the localization of $T$ in the point $x$.
This remark together with Lemma \ref{torsion} proves the Proposition.
\end{proof}

The above Proposition implies, that in order to describe general fiber
$F(\EE_\bullet,D)$ it is enough to have a description of the fibers
$F(\EE_\bullet,\gamma x)$, which we will call {\em simple fibers}.

\subsection{The stratification of a simple fiber}

We will need the following obvious Lemma.
\begin{lem}\label{coin}
Let $\EE$ be a vector bundle on $C$. Let $\EE'\subset\EE$ be a vector
subbundle, and let $E\subset\EE$ be a (necessarily locally free)
subsheaf. Then $E'=\EE'\cap E$ is a vector subbundle in $E$.

Moreover, the commutative square
$$
\begin{CD}
E'      @>>>    E       \\
@VVV        @VVV    \\
\EE'    @>>>    \EE
\end{CD}
$$
can be extended to the commutative diagram
$$
\begin{CD}
E'  @>>>    E   @>>>    E/E'    \\
@VVV        @VVV        @VVV    \\
\EE'    @>>>    \EE @>>>    \EE/\EE'\\
@VVV        @VVV        @VVV    \\
\EE'/E' @>>>    \EE/E   @>>> \dfrac{\EE/E}{\EE'/E'}\cong\dfrac{\EE/\EE'}{E/E'}
\end{CD}
$$
in which both the rows and the columns form the short exact sequences.
\end{lem}

The sheaf in the lower-right corner of the diagram
will be called {\em cointersection} of $E$ and $\EE'$ inside $\EE$
and denoted by $\coin\EE E{\EE'}$.

Let
$$
\gamma=\sum_{k=1}^{n-1}c_ki_k.
$$
For every $E_\bullet\in F(\EE_\bullet,\gamma x)$ we define
\begin{equation}
\mu_{pq}(E_\bullet)\ede l\left(\frac{\EE_q}{E_p\cap\EE_q}\right)
\qquad(1\le q\le p\le n-1),
\end{equation}
\begin{equation}\label{nu}
\nu_{pq}(E_\bullet)=\begin{cases}
\mu_{pq}(E_\bullet)-\mu_{p+1,q}(E_\bullet),&\text{if }1\le q\le p<n-1\\
\mu_{pq}(E_\bullet),&\text{if }1\le q\le p=n-1
\end{cases}
\end{equation}
\begin{equation}\label{ka}
\ka_{pq}(E_\bullet)=\begin{cases}
\nu_{pq}(E_\bullet)-\nu_{p,q-1}(E_\bullet),&\text{if }1<q\le p\le n-1\\
\nu_{pq}(E_\bullet),&\text{if }1=q\le p\le n-1
\end{cases}
\end{equation}
\begin{remark}
The transformations \refeq{nu} and \refeq{ka} are invertible, so the numbers
$\mu_{pq}$ can be uniquely reconstructed from $\nu_{pq}$ or $\ka_{pq}$.
Namely,
\begin{equation}\label{inv}
\ds\nu_{pq}=\sum_{r=1}^q\ka_{pr};\qquad
\ds\mu_{pq}=\sum_{s=p}^{n-1}\nu_{sq}=\sum_{r\le q\le p\le s}\ka_{sr}.
\end{equation}
\end{remark}

\begin{lem}
We have
\begin{equation}\label{pnu}
\nu_{pq}(E_\bullet)=l\left(\dfrac{\EE_q\cap E_{p+1}}{\EE_q\cap E_p}\right).
\end{equation}
\begin{equation}\label{pka}
\ka_{pq}(E_\bullet)=l\left(
\coin{\EE_q\cap E_{p+1}}{\EE_q\cap E_p}{\EE_{q-1}\cap E_{p+1}}\right).
\end{equation}
\end{lem}
\begin{proof}
The commutative diagram with exact rows
$$
\begin{CD}
0@>>> \EE_q\cap E_p     @>>> \EE_q @>>> \dfrac{\EE_q}{\EE_q\cap E_p}   @>>>0\\
@.  @VVV            @|      @VVV\\
0@>>> \EE_q\cap E_{p+1} @>>> \EE_q @>>> \dfrac{\EE_q}{\EE_q\cap E_{p+1}} @>>>0
\end{CD}
$$
implies \refeq{pnu}. In order to prove \refeq{pka} note that
$$
\EE_{q-1}\cap E_p=(\EE_q\cap E_p)\cap(\EE_{q-1}\cap E_{p+1})
$$
and apply Lemma \ref{coin} and \refeq{pnu}.
\end{proof}

\begin{cor}
Numbers  $\mu_{pq}$, $\nu_{pq}$ and $\ka_{pq}$
satisfy the following inequalities:
\begin{eqnarray}
0\le\ka_{pq}\\
0\le\nu_{p1}\le\nu_{p2}\le\dots\le\nu_{pp}\\
0\le\mu_{n-1,q}\le\mu_{n-2,q}\le\dots\le\mu_{qq}=c_q.\label{muin}
\end{eqnarray}
\end{cor}
\begin{proof}
See \refeq{pka},\refeq{ka},\refeq{nu} and compare the definition of $\mu_{qq}$
with \ref{torsion}.
\end{proof}

We will denote by $[p,q]$ the positive coroot
\begin{equation}\label{pq}
[p,q]\ede\sum_{k=q}^pi_k\in R^+
\end{equation}
\begin{lem}\label{part}
For any $E_\bullet\in F(\EE_\bullet,\gamma x)$ we have
$$
\sum_{1\le q\le p\le n-1}\ka_{pq}(E_\bullet)[p,q]=\gamma.
$$
\end{lem}
\begin{proof}
Applying \refeq{ka}, \refeq{pq} and \refeq{inv} we get
\begin{multline*}
\sum_{1\le q\le p\le n-1}\ka_{pq}[p,q]=
\sum_{1\le q\le p\le n-1}(\nu_{pq}-\nu_{p,q-1})[p,q]=\\=\!\!\!
\sum_{1\le q\le p\le n-1}\!\!\!\nu_{pq}([p,q]-[p,q+1])=\!\!\!
\sum_{1\le q\le p\le n-1}\!\!\!\nu_{pq}i_q=
\sum_{q=1}^{n-1}\left(\sum_{p=q}^{n-1}\nu_{pq}\right)i_q=
\sum_{q=1}^{n-1}\mu_{qq}i_q.
\end{multline*}
Now Lemma follows from \refeq{muin}.
\end{proof}

Let $\KK(\gamma)$ be the set of all partitions of $\gamma\in\NN[I]$ into the
sum of positive coroots: $\gamma=\suml_{s=1}^t\delta_s$, where
$\delta_s\in R^+$ (note that $\KK(\gamma)\ne\Gamma(\gamma)$).
In other words, since every positive coroot for $G=SL(V)$ is
equal to $[p,q]$ for some $p,q$,
$$
\KK(\gamma)=\{ (\ka_{pq})_{1\le q\le p\le n-1}\ |\
\ka_{pq}\ge0\quad\text{ and }\sum_{1\le q\le p\le n-1}\ka_{pq}[p,q]=\gamma\}.
$$
Let $\MMM(\gamma)$ denote the set of all collections $(\mu_{pq})$
which can be produced from some $(\ka_{pq})\in\KK(\gamma)$ as in \refeq{inv}.

The Lemma \ref{part} implies that for any
$E_\bullet\in F(\EE_\bullet,\gamma x)$ we have
$(\mu_{pq}(E_\bullet))\in\MMM(\gamma)$.
Define the stratum $\SS((\mu_{pq})_\pq{n-1},(\EE_k)_{k=1}^{n-1})$ as follows:
$$
\SS((\mu_{pq})_\pq{n-1},(\EE_k)_{k=1}^{n-1})=
\{ E_\bullet\in F(\EE_\bullet,\gamma x)\ |\ \mu_{pq}(E_\bullet)=\mu_{pq}\}.
$$
To unburden the notations in the cases when it is clear which flag
$\EE_\bullet$ is used we will write just $\SS_\mu$.
We have obviously
\begin{equation}\label{stratification}
F(\EE_\bullet,\gamma x)=\bigsqcup_{\mu\in\MMM(\gamma)}\SS_\mu.
\end{equation}

\begin{remark}\label{shortflag}
We will also use the similar
varieties $\SS((\mu_{pq})_\pq{N},(\EE_k)_{k=1}^N)$ that can
be defined in the same way for any short flag $(\EE_k)_{k=1}^N$
(that is the flag of subbundles $\EE_1\subset\dots\subset\EE_N\subset\cV$
with $\rank\EE_k=k$).
\end{remark}

\subsection{The strata $\SS_\mu$}

In order to study $\SS_\mu$ we will introduce some more varieties.

For every $1\le N\le n-1$, a short flag of subbundles $(\EE_k)_{k=1}^N$
(see Remark \ref{shortflag}) and a collection of numbers $(\nu_k)_{k=1}^N$
such that $0\le\nu_1\le\dots\le\nu_N$, we define the space
$\TT((\nu_k)_{k=1}^N,(\EE_k)_{k=1}^N)$ as follows:
\begin{equation}\label{TT}
\TT((\nu_k)_{k=1}^N,(\EE_k)_{k=1}^N)=\{E\subset\EE_N\ |\
\rank(E)=N\quad\text{ and }\quad
l{\left(\frac{\EE_k}{\EE_k\cap E}\right)}=\nu_k\}
\end{equation}

We define {\em pseudoaffine} spaces by induction in dimension.
First, the affine line $\AAA^1$ is a pseudoaffine space. Now
a space $A$ is called pseudoaffine if it admits a fibration
$A\to B$ with pseudoaffine fibers and pseudoaffine $B$.

\begin{thm}
The space $\TT((\nu_k)_{k=1}^N,\!(\EE_k)_{k=1}^N)$ is
pseudoaffine of dimension $\!\suml_{k=1}^{N-1}\!\!\nu_k$.
\end{thm}
\begin{proof}
We use induction in $N$.
The case $N=1$ is trivial. There is only one subsheaf $E$ in the
line bundle $\EE_1$ with $\leng(\EE_1/E)=\nu_1$, namely
$E=\mm_x^{\nu_1}\cdot\EE_1$. This means that $\TT(\nu_1,\EE_1)$
is a point and the base of induction follows.

If $N>1$ then consider the map
$$
\tau:\TT((\nu_k)_{k=1}^N,(\EE_k)_{k=1}^N)\to
\TT((\nu_k)_{k=1}^{N-1},(\EE_k)_{k=1}^{N-1}),
$$
which sends $E\in\TT((\nu_k)_{k=1}^N,(\EE_k)_{k=1}^N)$ to
$E'=E\cap\EE_{N-1}\in\TT((\nu_k)_{k=1}^{N-1},(\EE_k)_{k=1}^{N-1})$.
\begin{lem}\label{fibre}
Let $L=\mm_x^{\nu_N-\nu_{N-1}}\cdot\left(\dfrac{\EE_N}{\EE_{N-1}}\right)$.
For any $E'\in\TT((\nu_k)_{k=1}^{N-1},(\EE_k)_{k=1}^{N-1})$
there is an isomorphism
$$
\tau^{-1}(E')\cong\Hom(L,\EE_{N-1}/E')\cong
\AAA^{\leng(\EE_{N-1}/E')}=\AAA^{\nu_{N-1}}.
$$
\end{lem}
Thus, the space $\TT((\nu_k)_{k=1}^N,(\EE_k)_{k=1}^N)$ is
affine fibration over a pseudoaffine space, hence it is pseudoaffine and its
dimension is equal to
$$
\dim\left(\TT((\nu_k)_{k=1}^{N-1},(\EE_k)_{k=1}^{N-1})\right)+\nu_{N-1}=
\sum_{k=1}^{N-2}\nu_k+\nu_{N-1}=\sum_{k=1}^{N-1}\nu_k.
$$
The Theorem is proved.
\end{proof}
\begin{proofof}{Lemma \ref{fibre}}
Let $E\in\tau^{-1}(E')$. Since $E'=E\cap\EE_{N-1}$ we can apply Lemma
\ref{coin} which gives the following commutative diagram:
$$
\begin{CD}
E'      @>>>    E   @>>>    L       \\
@VVV            @VVV        @VVV        \\
\EE_{N-1}   @>j>>   \EE_N   @>\psi>>    \EE_N/\EE_{N-1} \\
@VVV            @VVV        @VVV        \\
T_{N-1}     @>>>    T_N @>>>    T_N/T_{N-1}
\end{CD}
$$
(Note that since $\EE_N/\EE_{N-1}$ is a line bundle and
$\leng(T_N/T_{N-1})=\leng(T_N)-\leng(T_{N-1})=\nu_N-\nu_{N-1}$ the kernel
of natural map $\EE_N/\EE_{N-1}\to T_N/T_{N-1}$ is isomorphic
to $L$.)
Let $\ti\EE_N=\psi^{-1}(L)$. Then $E$ is contained in $\ti\EE_N$ and
we have the following commutative diagram:
$$
\begin{CD}
E'      @>>>    E       @>>>        L   \\
@VVV            @VVV                @|  \\
\EE_{N-1}       @>j>>   \ti\EE_N        @>\psi>>        L       \\
@VVV            @V\eps VV               \\
T_{N-1}     @=  T_{N-1}
\end{CD}
$$
This means that the points of $\tau^{-1}(E')$ are in one-to-one correspondence
with maps $\eps:\ti\EE_N\to T_{N-1}$ such that $\eps\cdot j$ is equal
to the canonical projection from $\EE_{N-1}$ to $T_{N-1}$.
Applying the functor $\Hom(\bullet,T_{N-1})$ to the middle row of the
above diagram we get an exact sequence:
$$
0\to\Hom(L,T_{N-1})\to\Hom(\ti\EE_N,T_{N-1})@>j^*>>\Hom(\EE_{N-1},T_{N-1})\to
\Ext^1(L,T_{N-1}).
$$
The last term in this sequence is zero because $L$ is locally free
and $T_{N-1}$ is torsion.
This means that the space of maps $\eps$ which we need to describe
is a torsor over the group $\Hom(L,T_{N-1})$. Hence this space can
be identified with the group. Thus, we have proved that
$\tau^{-1}(E')\cong\Hom(L,T_{N-1})$ is an affine space.

Now,
$$
\dim(\tau^{-1}(E'))=\dim\Hom(L,T_{N-1})=\dim\HH^0(T_{N-1})=\leng(T_{N-1})=
\nu_{N-1}.
$$
The Lemma is proved.
\end{proofof}

\begin{thm}\label{saff}
The space $\SS((\mu_{pq})_\pq{N},(\EE_k)_{k=1}^N)$
is a pseudoaffine space of dimension $\mu_{21}+\mu_{32}+\dots+\mu_{N,N-1}$.
\end{thm}
\begin{proof}
We use induction in $N$.
If $N=1$ then $\SS_\mu$ is a point and the base of induction follows.

If $N>1$ consider the map
$$
\sigma:\SS((\mu_{pq})_\pq{N},(\EE_k)_{k=1}^N)\to
\TT((\mu_{N,k})_{k=1}^N,(\EE_k)_{k=1}^N),
$$
which sends $(E_k)_{k=1}^N$ to $E_N\subset\EE_N$.
\begin{lem}\label{sfibre}
Let $E\in\TT((\nu_{N,k})_{k=1}^N,(\EE_k)_{k=1}^N)$.
Consider $\ti\EE_k=\EE_k\cap E\quad(1\le k\le N-1)$ and
set $\ti\mu_{pq}=\mu_{pq}-\mu_{Nq}\quad(\pq{N-1})$.
Then $(\ti\EE_k)_{k=1}^{N-1}$ is a short flag of subbundles and
for any $E\in\TT((\mu_{N,k})_{k=1}^N,(\EE_k)_{k=1}^N$
\begin{equation}
\sigma^{-1}(E)\cong\SS((\ti\mu_{pq})_{\pq{N-1}}),(\ti\EE_k)_{k=1}^{N-1}).
\end{equation}
\end{lem}
Thus $\SS((\mu_{pq})_\pq{N},(\EE_k)_{k=1}^N)$ is a fiber space
with pseudoaffine base and fiber, therefore it is pseudoaffine.

Now, the calculation of the dimension
\begin{multline*}
\dim\left(\SS((\mu_{pq})_\pq{N},(\EE_k)_{k=1}^N)\right)=
\sum_{k=1}^{N-1}\mu_{N,k}+\sum_{k=1}^{N-2}\ti\mu_{k+1,k}=\\=
\sum_{k=1}^{N-1}\mu_{N,k}+\sum_{k=1}^{N-2}(\mu_{k+1,k}-\mu_{N,k})=
\mu_{N,N-1}+\sum_{k=1}^{N-2}\mu_{k+1,k}=
\sum_{k=1}^{N-1}\mu_{k+1,k},
\end{multline*}
finishes the proof of the Theorem.
\end{proof}

\begin{proofof}{Lemma \ref{sfibre}}
Assume that $(E_k)_{k=1}^N\in\SS_\mu$ and $E_N=E$. The commutative
diagram
$$
\begin{CD}
0 @>>>  \EE_q\cap E_p   @>>> \EE_q @>>> \dfrac{\EE_q}{\EE_q\cap E_p} @>>> 0\\
@.      @VVV         @|     @VVV                  @.\\
0 @>>>  \EE_q\cap E @>>> \EE_q @>>> \dfrac{\EE_q}{\EE_q\cap E}   @>>> 0
\end{CD}
$$
implies that
\begin{equation}\label{lengths}
l\left(\dfrac{\EE_q\cap E}{\EE_q\cap E_p}\right)=
l\left(\dfrac{\EE_q}{\EE_q\cap E_p}\right)-
l\left(\dfrac{\EE_q}{\EE_q\cap E}\right)=
\mu_{pq}-\mu_{Nq}=\ti\mu_{pq},
\end{equation}
hence $(E_k)_{k=1}^{N-1}\in\SS_{\ti\mu}$.

Vice versa, if $(E_k)_{k=1}^{N-1}\in\SS_{\ti\mu}$ then
the above commutative diagram along with \refeq{lengths} implies that
$(E_k)_{k=1}^N\in\SS_\mu$, where we have put $E_N=E$.
\end{proofof}

\subsection{The cohomology of the simple fiber}
Now we will compute the dimension of the strata $\SS_\mu$ in
terms of the partition $\ka$.

\begin{defn}
A space $\cX$ is called {\em cellular} if it
admits a stratification with pseudoaffine strata.
\end{defn}

Suppose $\cX=\bigsqcup\limits_{\xi\in\Xi} S_\xi$ is a
pseudoaffine stratification of a cellular space $\cX$.
For a positive integer $j$ we define
$\chi(j)\ede\{\xi\in\Xi\ |\ \dim S_\xi=j\}.$

\begin{lem}\label{hodge}
The Hodge structure $\HH^\bullet(\cX,\QQ)$ is a direct sum of Tate structures,
and $\QQ(j)$ appears with multiplicity $\chi(j)$.
In other words,
$$
\HH^\bullet(\cX,\QQ)=\oplus_{j\in\NN}\QQ(j)^{\chi(j)}.
$$
\end{lem}
\begin{proof}
Evident.
\end{proof}

Given a Tate structure $\CH=\oplus_{j\in\NN}\QQ(j)^{\chi(j)}$ we
consider a {\em generating function}
$$
P(\CH,t)=\sum_{j\in\NN}\chi(j)t^j\in\NN[t].
$$

For $\ka\in\KK(\gamma)$ we define $K(\ka)\ede\ds\suml_{\pq{n-1}}\ka_{pq}$
as the number of summands in the partition $\ka$.
For $\gamma\in\NN[I]$ the following $q$-analog of the Kostant's
partition function was was introduced in \cite{Lusztig}:
\begin{equation}\label{cke}
\cK_\gamma(t)=t^{|\gamma|}\sum_{\ka\in\KK(\gamma)}t^{-K(\ka)}.
\end{equation}

\begin{lem}\label{dims}
Let $\ka\in\KK(\gamma)$ and $\mu\in\MMM(\gamma)$ be defined
as in \refeq{inv}. Then
\begin{equation}
\dim\SS_\mu=\sum_{k=1}^{n-2}\mu_{k+1,k}=|\gamma|-K(\ka),
\end{equation}
\end{lem}
\begin{proof}
Applying \refeq{inv} we get
\begin{multline*}
\sum_{k=1}^{n-2}\mu_{k+1,k}=
\sum_{k=1}^{n-2}
\left(\sum\begin{Sb}1\le q\le k\\k+1\le p\le n-1\end{Sb}\ka_{pq}\right)=
\sum_{\pq{n-1}}(p-q)\ka_{pq}=\\=
\sum_{\pq{n-1}}(|[p,q]|-1)\ka_{pq}=
|\gamma|-\sum_{\pq{n-1}}\ka_{pq}=|\gamma|-K(\ka).
\end{multline*}
{}\end{proof}

\begin{cor}\label{ck}
For any $\gamma\in\NN[I]$, $x\in C$, the simple fiber
$F(\EE_\bullet,\gamma x)$ is a cellular space, and the generating
function of its cohomology is equal to the Lusztig--Kostant polynomial
$$
P(\HH^\bullet(F(\EE_\bullet,\gamma x)),t)=\cK_\gamma(t).
$$
\end{cor}
\begin{proof}
Apply \refeq{stratification}, \ref{saff}, \ref{hodge} and \ref{dims}.
\end{proof}

\begin{cor}\label{cck}
Let $D=\suml_{r=1}^m\gamma_rx_r\in C^{\alpha-\beta}_\Gamma$.
The fiber $F(\EE_\bullet,D)$ is a cellular space and
\begin{equation}\label{ccke}
P(\HH^\bullet(F(\EE_\bullet,D)),t)=\cK_\Gamma
\ede\prod_{r=1}^m\cK_{\gamma_r}(t).
\end{equation}
\end{cor}
\begin{proof}
Apply \ref{gfibre}, \ref{hodge} and \ref{ck}.
\end{proof}

\begin{lem}\label{est}
Let $D=\suml_{r=1}^m\gamma_r x_r$. We have
$$
\dim F(\EE_\bullet,D)\le\left|\sum_{r=1}^m\gamma_r\right|-m.
$$
\end{lem}
\begin{proof}
Note that for any $\ka\in\KK(\gamma_r)$ we have $K(\ka)\ge1$, hence
$\deg\cK_{\gamma_r}\le|\gamma_r|-1$. Now, the Lemma follows from \ref{cck}.
\end{proof}

\begin{proofof}{Main Theorem}
Consider the stratum $\DD\beta\Gamma$ of $\MMD\alpha$.
Its dimension is $2|\beta|+\dim\FF+m$ and codimension is
$2|\alpha-\beta|-m$. The Lemma \ref{est} implies that the dimension
of the fiber of $\pi$ over the stratum $\DD\beta\Gamma$ is less than or equal
to $|\alpha-\beta|-m$, which is strictly less then the half
codimension of the stratum.
\end{proofof}

\subsection{Applications}
Let $\UQ$ denote the smooth constant Hodge  irreducible module on $\MML\alpha$
(as a constructible complex it lives in cohomological degree
$-2|\alpha|-\dim\FF$). Let $IC$ denote the minimal extension
of a smooth constant irreducible Hodge module from $\MM\alpha$ to
$\MMD\alpha$. It is well known that the smallness of $\pi$ implies
the following corollary.
\begin{cor}
$$
IC=\pi_*\UQ.
$$
\end{cor}

Now we can compute the stalks of $IC$ as
cohomology of fibers of $\pi$:
for $\vphi\in\MMD\alpha$ we have
$$
IC_{(\vphi)}^\bullet=\HH^\bullet(\pi^{-1}(\vphi),\UQ)
$$
as graded Hodge structures.

\begin{cor}[Parity vanishing]
$$
IC_{(\vphi)}^{j}=0\quad\text{if $j-\dim\FF$ is odd.}
$$
\end{cor}
\begin{proof}
Use \ref{cck}.
\end{proof}

\begin{cor}
For $\vphi\in\DD\beta\Gamma$ we have
$$
IC_{(\vphi)}^{-2|\alpha|-\dim\FF+2j}=\QQ(j)^{\fk_\Gamma(j)},
$$
where $\fk_\Gamma(j)$ is the coefficient of $t^j$ in $\KK_\Gamma(t)$.
\end{cor}


\begin{thebibliography}{XXXX}

\bibitem[Dri]{Drinfeld} Drinfeld V.,
Two-dimensional $\ell$-adic representations of the fundamental
group of a curve over a finite field and automorphic forms on $GL(2)$,
American Journal of Mathematics, {\bf 105}, 1983, pp.~85--114.

\bibitem[Giv]{Givental} Givental A.,
Equivariant Gromow--Witten Invariants,
preprint, alg-geom/9603021.

\bibitem[Kon]{Kontsevich} Kontsevich M.,
Enumeration of rational curves via torus actions,
In: The moduli space of curves,
R.~Dijkgraaf, C.~Faber, G.~van der Geer (Eds.),
Progress in Math., vol.~129, Birkh\"auser, 1995, pp.~335--368.

\bibitem[Lau]{Laumon} Laumon G.,
Faisceaux Automorphes Li\'es aux S\'eries d'Eisenstein,
In: Automorphic Forms, Shimura Varietes, and $L$-functions,
Perspect.\ Math., vol.~10, Academic Press, Boston, MA, 1990, pp.~227--281

\bibitem[Lus]{Lusztig} Lusztig G.,
Singularities, character formulas and a $q$-analog of weight multiplicities,
Ast\'erisque {\bf 101--102}, 1983, pp.~208--229.

\end{thebibliography}
\enddocument